\input amstex
\documentstyle{amsppt}
%
\catcode`@=11
\redefine\output@{%
  \def\break{\penalty-\@M}\let\par\endgraf
  \ifodd\pageno\global\hoffset=105pt\else\global\hoffset=8pt\fi  
  \shipout\vbox{%
    \ifplain@
      \let\makeheadline\relax \let\makefootline\relax
    \else
      \iffirstpage@ \global\firstpage@false
        \let\rightheadline\frheadline
        \let\leftheadline\flheadline
      \else
        \ifrunheads@ 
        \else \let\makeheadline\relax
        \fi
      \fi
    \fi
    \makeheadline \pagebody \makefootline}%
  \advancepageno \ifnum\outputpenalty>-\@MM\else\dosupereject\fi
}
\catcode`\@=\active
\nopagenumbers
\def\negskp{\hskip -2pt}
\def\MatGrU{\operatorname{U}}
\def\MatGrSU{\operatorname{SU}}
\def\MatAlgSU{\operatorname{su}}
\def\Alpha{\operatorname{A}}
\def\Re{\operatorname{Re}}
\def\vtrule{\vrule height 12pt depth 6pt}
\def\chirk{\special{em:point 1}\kern 1.2pt\raise 0.6pt
  \hbox to 0pt{\special{em:point 2}\hss}\kern -1.2pt
  \special{em:line 1,2,0.3pt}\ignorespaces}
\def\Chirk{\special{em:point 1}\kern 1.5pt\raise 0.6pt
  \hbox to 0pt{\special{em:point 2}\hss}\kern -1.5pt
  \special{em:line 1,2,0.3pt}\ignorespaces}
\def\Hirk{\kern 0pt\special{em:point 1}\kern 4pt\special{em:point 2}
  \kern -3pt\special{em:line 1,2,0.3pt}\ignorespaces}
\accentedsymbol\uuud{d\unskip\kern -3.8pt\raise 1pt\hbox to
  0pt{\chirk\hss}\kern 0.2pt\raise 1.9pt\hbox to 0pt{\chirk\hss}
  \kern 0.2pt\raise 2.8pt\hbox to 0pt{\chirk\hss}\kern 3.5pt}
\accentedsymbol\bolduuud{\bold d\unskip\kern -4pt\raise 1.1pt\hbox to
  0pt{\chirk\hss}\kern 0pt\raise 1.9pt\hbox to 0pt{\chirk\hss}
  \kern 0pt\raise 2.5pt\hbox to 0pt{\chirk\hss}\kern 4pt}
\accentedsymbol\uud{d\unskip\kern -3.72pt\raise 1.5pt\hbox to
  0pt{\chirk\hss}\kern 0.26pt\raise 2.4pt\hbox to 0pt{\chirk\hss}
  \kern 3.5pt}
\accentedsymbol\bolduud{\bold d\unskip\kern -4pt\raise 1.4pt\hbox to
  0pt{\chirk\hss}\kern 0pt\raise 2.3pt\hbox to 0pt{\chirk\hss}
  \kern 4pt}
\accentedsymbol\uuuD{D\unskip\kern -4.6pt\raise 2pt\hbox to
  0pt{\Chirk\hss}\kern 0.2pt\raise 3pt\hbox to 0pt{\Chirk\hss}
  \kern 0.2pt\raise 4.1pt\hbox to 0pt{\Chirk\hss}\kern 4.5pt}
\accentedsymbol\uuD{D\unskip\kern -4.5pt\raise 2.4pt\hbox to
  0pt{\Chirk\hss}\kern 0.2pt\raise 3.4pt\hbox to 0pt{\Chirk\hss}
  \kern 4.4pt}
\accentedsymbol\uD{D\unskip\kern -4.4pt\raise 3pt\hbox to
  0pt{\Chirk\hss}\kern 4.4pt}
\accentedsymbol\bolduuuD{\bold D\unskip\kern -4.6pt\raise 2pt\hbox to
  0pt{\Chirk\hss}\kern 0pt\raise 3pt\hbox to 0pt{\Chirk\hss}
  \kern 0pt\raise 4.1pt\hbox to 0pt{\Chirk\hss}\kern 4.5pt}
\accentedsymbol\bolduuD{\bold D\unskip\kern -4.5pt\raise 2.4pt\hbox to
  0pt{\Chirk\hss}\kern 0pt\raise 3.4pt\hbox to 0pt{\Chirk\hss}
  \kern 4.4pt}
\accentedsymbol\bolduD{\bold D\unskip\kern -4.5pt\raise 3pt\hbox to
  0pt{\Chirk\hss}\kern 4.4pt}
\accentedsymbol\uuuU{U\unskip\kern -5.1pt\raise 2pt\hbox to
  0pt{\Chirk\hss}\kern 0.2pt\raise 3pt\hbox to 0pt{\Chirk\hss}
  \kern 0.2pt\raise 4.1pt\hbox to 0pt{\Chirk\hss}\kern 4.9pt}
\accentedsymbol\uuU{U\unskip\kern -5.0pt\raise 2.4pt\hbox to
  0pt{\Chirk\hss}\kern 0.2pt\raise 3.4pt\hbox to 0pt{\Chirk\hss}
  \kern 4.9pt}
\accentedsymbol\uU{U\unskip\kern -4.9pt\raise 3pt\hbox to
  0pt{\Chirk\hss}\kern 4.8pt}
\accentedsymbol\bolduPsi{\boldsymbol\Psi\unskip\kern -6.5pt\raise 3pt
  \hbox to 0pt{\Hirk\hss}\kern 6.6pt}
\accentedsymbol\bolduuPsi{\boldsymbol\Psi\unskip\kern -6.5pt\raise 3pt
  \hbox to 0pt{\Hirk\hss}\unskip\kern 0pt\raise 3.8pt
  \hbox to 0pt{\Hirk\hss}\kern 6.6pt}
\accentedsymbol\bolduuuPsi{\boldsymbol\Psi\unskip\kern -6.5pt\raise 3pt
  \hbox to 0pt{\Hirk\hss}\unskip\kern 0pt\raise 3.8pt
  \hbox to 0pt{\Hirk\hss}\unskip\kern 0pt\raise 4.6pt
  \hbox to 0pt{\Hirk\hss}\kern 6.6pt}
\accentedsymbol\uA{\operatorname{A}\unskip\kern -6.2pt\raise 2.8pt\hbox to
  0pt{\vbox{\hrule width 2pt height 0.3pt}\hss}\kern 4.5pt}
\accentedsymbol\uuA{\operatorname{A}\unskip\kern -6.2pt\raise 2.8pt\hbox 
  to 0pt{\vbox{\hrule width 2pt height 0.3pt}\hss}\unskip\kern 0.2pt\raise
  3.4pt\hbox to 0pt{\vbox{\hrule width 1.6pt height 0.3pt}\hss}\kern 4.3pt}
\accentedsymbol\uuuA{\operatorname{A}\unskip\kern -6.2pt\raise 2.8pt\hbox 
  to 0pt{\vbox{\hrule width 2pt height 0.3pt}\hss}\unskip\kern 0.2pt\raise
  3.4pt\hbox to 0pt{\vbox{\hrule width 1.6pt height 0.3pt}\hss}\unskip
  \kern 0.2pt\raise 4.0pt\hbox to 0pt{\vbox{\hrule width 1.0pt height
  0.3pt}\hss}\kern 4.3pt}
\def\blue#1{#1}
\catcode`#=11\def\diez{#}\catcode`#=6
\catcode`_=11\def\podcherkivanie{_}\catcode`_=8
\def\mycite#1{\cite{\blue{#1}}\immediate\special{ps:
     ShrHPSdict begin /ShrBORDERthickness 0 def}}
\def\myciterange#1#2#3#4{\cite{\blue{#2#3#4}}\immediate\special{ps:
     ShrHPSdict begin /ShrBORDERthickness 0 def}}
\def\mytag#1{%
    \tag#1}
\def\mythetag#1{\thetag{\blue{#1}}\immediate\special{ps:
     ShrHPSdict begin /ShrBORDERthickness 0 def}}
\def\myrefno#1{\no#1}
\def\myhref#1#2{\blue{#2}\immediate\special{ps:
     ShrHPSdict begin /ShrBORDERthickness 0 def}}
\def\myEarXivlink{\myhref{http://arXiv.org}{http:/\negskp/arXiv.org}}
\def\myGeoCities{\myhref{http://www.geocities.com}{GeoCities}}
\def\mytheorem#1{\csname proclaim\endcsname{Theorem #1}}

\def\mylemma#1{\csname proclaim\endcsname{Lemma #1}}

\def\mycorollary#1{\csname proclaim\endcsname{Corollary #1}}

\def\mydefinition#1{\definition{Definition #1}}
\def\mythedefinition#1{\blue{#1}\immediate\special{ps:
     ShrHPSdict begin /ShrBORDERthickness 0 def}}

\pagewidth{360pt}
\pageheight{606pt}
\topmatter
\title
A note on connections of the Standard Model 
in a gravitation field.
\endtitle
\author
R.~A.~Sharipov
\endauthor
\address 5 Rabochaya street, 450003 Ufa, Russia\newline
\vphantom{a}\kern 12pt Cell Phone: +7-(917)-476-93-48
\endaddress
\email \vtop to 30pt{\hsize=280pt\noindent
\myhref{mailto:r-sharipov\@mail.ru}
{r-sharipov\@mail.ru}\newline
\myhref{mailto:ra\podcherkivanie sharipov\@lycos.com}{ra\_\hskip 1pt
sharipov\@lycos.com}\newline
\myhref{mailto:R\podcherkivanie Sharipov\@ic.bashedu.ru}
{R\_\hskip 1pt Sharipov\@ic.bashedu.ru}\vss}
\endemail
\urladdr
\vtop to 20pt{\hsize=280pt\noindent
\myhref{http://www.geocities.com/r-sharipov}
{http:/\negskp/www.geocities.com/r-sharipov}\newline
\myhref{http://www.freetextbooks.boom.ru/index.html}
{http:/\negskp/www.freetextbooks.boom.ru/index.html}\vss}
\endurladdr
\abstract
    The Standard Model of the theory of elementary particles is based on 
the $\MatGrU(1)\times\MatGrSU(2)\times\MatGrSU(3)$ symmetry. In the presence
of a gravitation field, i\.\,e\. in a non-flat space-time manifold, this
symmetry is implemented through three special vector bundles. Connections
associated with these vector bundles are studied in this paper. In the
Standard Model they are interpreted as gauge fields.
\endabstract
\endtopmatter
\loadbold
\loadeufb
\TagsOnRight
\document

\rightheadtext{On connections of the Standard Model \dots}
\head
1. The $\MatGrU$ \!and $\MatGrSU$-bundles and their basic fields. 
\endhead
    In its canonical form the Standard Model describes elementary 
particles in the flat Minkowski space-time (see \myciterange{1}1{--\,}6). 
When passing to a non-flat space-time $M$ we introduce three special 
vector bundles $\uU\!M$, $S\uuU\!M$, and $S\uuuU\!M$, each equipped with
its own basic tensorial fields (see \mycite{7}). They are listed in the
following table:
$$
\vcenter{\hsize 10cm
\offinterlineskip\settabs\+\indent
\vtrule
\hskip 1.2cm &\vtrule 
\hskip 4.2cm &\vtrule 
\hskip 1.7cm &\vtrule 
\hskip 2.7cm &\vtrule 
\cr\hrule 
\+\vtrule
\hfill\,Symbol\hfill&\vtrule
\hfill Name\hfill &\vtrule
\hfill Bundle\hfill &\vtrule
\hfill Tensorial\hfill &\vtrule\cr
\vskip -0.2cm
\+\vtrule
\hfill &\vtrule
\hfill \hfill&\vtrule
\hfill \hfill&\vtrule
\hfill type\hfill&\vtrule\cr\hrule
\+\vtrule
\hfill $\bolduD$\hfill&\vtrule
\hfill Hermitian metric tensor\hfill&\vtrule
\hfill \uU\!M\hfill&\vtrule
\hfill $(0,1|0,1)$\hfill&\vtrule\cr\hrule
\+\vtrule
\hfill$\bolduuD$\hfill&\vtrule
\hfill Hermitian metric tensor\hfill&\vtrule
\hfill S\uuU\!M\hfill&\vtrule
\hfill $(0,1|0,1)$\hfill&\vtrule\cr
\+\vrule height 0.4pt depth 0pt width 1.63cm
&\vrule height 0.4pt depth 0pt width 4.23cm
&\hskip 1.7cm
&\vrule height 0.4pt depth 0pt width 2.73cm\cr
\+\vtrule
\hfill$\bolduud$\hfill&\vtrule
\hfill Skew-symmetric\hfill&\vtrule
\hfill \hfill&\vtrule
\hfill $(0,2|0,0)$\hfill&\vtrule\cr
\+\vtrule
\hfill &\vtrule\hfill metric tensor\hfill&\vtrule
\hfill &\vtrule\hfill &\vtrule\cr\hrule
\+\vtrule
\hfill$\bolduuuD$\hfill&\vtrule
\hfill Hermitian metric tensor\hfill&\vtrule
\hfill $S\uuuU\!M$\hfill&\vtrule
\hfill $(0,1|0,1)$\hfill&\vtrule\cr
\+\vrule height 0.4pt depth 0pt width 1.63cm
&\vrule height 0.4pt depth 0pt width 4.23cm
&\hskip 1.7cm
&\vrule height 0.4pt depth 0pt width 2.73cm\cr
\+\vtrule
\hfill$\bolduuud$\hfill&\vtrule
\hfill Completely\hfill&\vtrule
\hfill \hfill&\vtrule
\hfill $(0,3|0,0)$\hfill&\vtrule\cr
\+\vtrule
\hfill &\vtrule\hfill skew-symmetric tensor\hfill&\vtrule
\hfill &\vtrule\hfill &\vtrule\cr\hrule
}\quad
\mytag{1.1}
$$
The bundles $\uU\!M$, $S\uuU\!M$, and $S\uuuU\!M$ are complex
bundles over a real base. For this reason they possess the involution
of complex conjugation $\tau$:
$$
\gather
\hskip -2em
\CD
@>\tau>>\\
\vspace{-4ex}
\uU^\varepsilon_\eta\bar{\uU}^{\kern 0.2pt\lower 2.1pt\hbox{$\ssize
\sigma$}}_\zeta M@.
\uU^\sigma_\zeta\bar{\uU}^{\kern 0.2pt\lower 2.1pt
\hbox{$\ssize\varepsilon$}}_\eta M,\\
\vspace{-4.2ex}
@<<< 
\endCD
\mytag{1.2}\\
\hskip -2em
\CD
@>\tau>>\\
\vspace{-4ex}
S\uuU^\pi_\rho\overline{S\uuU}^{\kern 0.2pt\lower 2.1pt
\hbox{$\ssize\omega$}}_\mu M@.
S\uuU^\omega_\mu\overline{S\uuU}^{\kern 0.2pt\lower 2.1pt
\hbox{$\ssize\pi$}}_\rho M,\\
\vspace{-4.2ex}
@<<< 
\endCD
\mytag{1.3}\\
\hskip -2em
\CD
@>\tau>>\\
\vspace{-4ex}
S\uuuU^\varkappa_\lambda\overline{S\uuuU}^{\kern 0.2pt\lower 2.1pt
\hbox{$\ssize\varsigma$}}_\chi M@.
S\uuuU^\varsigma_\chi\overline{S\uuuU}^{\kern 0.2pt\lower 2.1pt
\hbox{$\ssize\varkappa$}}_\lambda M.\\
\vspace{-4.2ex}
@<<< 
\endCD
\mytag{1.4}
\endgather
$$
The tensor bundles $\uU^\varepsilon_\eta\bar{\uU}^{\kern 0.2pt
\lower 2.1pt\hbox{$\ssize\sigma$}}_\zeta M$, $S\uuU^\pi_\rho
\overline{S\uuU}^{\kern 0.2pt\lower 2.1pt\hbox{$\ssize
\omega$}}_\mu M$, and $S\uuuU^\varkappa_\lambda\overline{S
\uuuU}^{\kern 0.2pt\lower 2.1pt\hbox{$\ssize\varsigma$}}_\chi M$
from \mythetag{1.2}, \mythetag{1.3}, and \mythetag{1.4} are 
defined in \mycite{7} (see formulas \thetag{2.5}, \thetag{2.11}, 
and \thetag{3.4} there).\par
     The Hermitian metric tensors $\bolduD$, $\bolduuD$, and 
$\bolduuuD$ in the table \mythetag{1.1} are introduced through
the corresponding Hermitian forms:
$$
\gather
\hskip -2em
\uD(\bold X,\bold Y)=\uD_{11}\,\overline{X^1}\,Y^{\kern 0.2pt
\lower 1.2pt\hbox{$\ssize 1$}}=\overline{\uD(\bold Y,\bold X)},
\mytag{1.5}\\
\vspace{1ex}
\hskip -2em
\uuD(\bold X,\bold Y)=\sum^2_{i=1}\sum^2_{\bar j=1}
\uuD_{i\bar j}\,\overline{X^{\bar j}}\,Y^{\kern 0.2pt
\lower 1.2pt\hbox{$\ssize i$}}=\overline{\uuD(\bold Y,\bold X)},
\mytag{1.6}\\
\hskip -2em
\uuuD(\bold X,\bold Y)=\sum^3_{i=1}\sum^3_{\bar j=1}
\uuuD_{i\bar j}\,\overline{X^{\bar j}}\,Y^{\kern 0.2pt
\lower 1.2pt\hbox{$\ssize i$}}=\overline{\uuuD(\bold Y,\bold X)}.
\mytag{1.7}
\endgather
$$
The Hermitian forms \mythetag{1.5}, \mythetag{1.6}, and \mythetag{1.7} 
are positive, i\.\,e\. their signatures are $(+)$, $(+,+)$, and $(+,+,+)$
respectively. The components of tensor fields $\uD_{11}$, $\uuD_{i\bar j}$,
and $\uuuD_{i\bar j}$ in the above formulas are frame-relative.
\mydefinition{1.1}{\it A frame\/} of a $q$-dimensional vector bundle 
over the base manifold $M$ is an ordered set of $q$ smooth sections of 
this bundle linearly independent at each point $p$ of some open domain 
$U\subset M$.
\enddefinition
    Let's denote by $(U,\bolduPsi_1)$, $(U,\,\bolduuPsi_1,\,\bolduuPsi_2)$,
$(U,\,\bolduuuPsi_1,\,\bolduuuPsi_2,\,\bolduuuPsi_3)$ some frames of the
bundles $\uU\!M$, $S\uuU\!M$, and $S\uuuU\!M$ respectively. Apart from 
these three frames, we choose some frame $(U,\,\boldsymbol\Upsilon_0,\,
\boldsymbol\Upsilon_1,\,\boldsymbol\Upsilon_2,\,\boldsymbol\Upsilon_3)$
of the tangent bundle $TM$.
\mydefinition{1.2} A frame $(U,\bolduPsi_1)$ of the electro-weak bundle
$\uU\!M$ is called an {\it orthonormal frame} if $\uD_{11}=1$ in this 
frame.
\enddefinition
\mydefinition{1.3} A frame $(U,\,\bolduuPsi_1,\,\bolduuPsi_2)$ of the
bundle $S\uuU\!M$ is called an {\it orthonormal frame} if the tensor 
fields $\bolduuD$ and $\bolduud$ are given by the following matrices
in this frame:
$$
\xalignat 2
&\uuD_{i\bar j}=\Vmatrix 1 & 0\\0 & 1\endVmatrix,
&&\uud_{ij}=\Vmatrix 0 & 1\\-1 & 0\endVmatrix.
\endxalignat
$$
\enddefinition
\mydefinition{1.4} A frame $(U,\,\bolduuuPsi_1,\,\bolduuuPsi_2,
\,\bolduuuPsi_3)$ of the bundle $S\uuuU\!M$ is called an {\it orthonormal
frame} if $\uuud_{123}=1$ and $\bolduuuD$ is given by the unit 
matrix in this frame:
$$
\uuuD_{i\bar j}=\Vmatrix 1 & 0 & 0\\
\vspace{0.5ex}
0 & 1 & 0\\
\vspace{0.5ex}
0 & 0 & 1\endVmatrix.
$$
\enddefinition
Once the component $\uuud_{123}$ is fixed by the condition 
$\uuud_{123}=1$, other components of the tensor $\bolduuud$ 
are determined by the skew-symmetry condition:
$$
\xalignat 3
&\uuud_{ijk}=-\uuud_{ji\kern 0.4pt k},
&&\uuud_{ijk}=-\uuud_{i\kern 0.4pt kj},
&&\uuud_{ijk}=-\uuud_{kj\kern 0.4pt i}.
\qquad
\mytag{1.8}
\endxalignat
$$
Orthonormal frames of the bundle $\uU\!M$ do always exist. As for the
other two bundles $S\uuU\!M$ and $S\uuuU\!M$, orthonormal frames for
them do exist provided some concordance conditions are fulfilled (see
theorems~2.1 and 3.1 in \mycite{7}). In the case of the electro-weak 
bundle $S\uuU\!M$ the concordance condition of $\bolduuD$ and $\bolduud$ 
is written as
$$
\hskip -2em
\sum^2_{i=1}\sum^2_{j=1}\uud^{\kern 0.4pt ij}
\,\uuD_{i\bar i}\,\uuD_{j\bar j}
=\overline{\uud_{\kern 0.4pt\bar i\kern 0.2pt\bar j}}.
\mytag{1.9}
$$
By $\uud^{\kern 0.4pt ij}$ in \mythetag{1.9} we denote the components 
of the matrix inverse to the matrix $\uud_{ij}$. They define the inverse
skew-symmetric metric tensor in $S\uuU\!M$. By tradition this tensor is
denoted by the same symbol $\bolduud$ as the initial one.\par
     In the case of the color bundle $S\uuuU\!M$ the concordance condition
for the tensor fields $\bolduuuD$ and $\bolduuud$ is given by the following
equality:
$$
\hskip -2em
\sum^3_{i=1}\sum^3_{j=1}\sum^3_{k=1}\uuud^{\kern 0.4pt ijk}
\,\uuuD_{i\bar i}\,\uuuD_{j\bar j}\,\uuuD_{k\bar k}
=\overline{\uuud_{\kern 0.4pt\bar i\kern 0.2pt\bar j\bar k}}.
\mytag{1.10}
$$
The equality \mythetag{1.10} is a three-dimensional analog of the equality
\mythetag{1.9}. By $\uuud^{\kern 0.4pt ijk}$ in \mythetag{1.10} we denote
the components of the completely skew-symmetric tensor of the type $(3,0|
0,0)$ inverse to the tensor $\bolduuud$. This tensor is denoted by the same
symbol $\bolduuud$ as the initial one. Due to the skew-symmetry conditions
$$
\xalignat 3
&\uuud^{\kern 0.4pt ijk}=-\uuud^{ji\kern 0.4pt k},
&&\uuud^{\kern 0.4pt ijk}=-\uuud^{\kern 0.4pt i\kern 0.4pt kj},
&&\uuud^{\kern 0.4pt ijk}=-\uuud^{\kern 0.4pt kj\kern 0.4pt i}
\endxalignat
$$
similar to \mythetag{1.8} the inverse tensor $\bolduuud$ is fixed by
the equality $\uuud^{\kern 0.4pt 123}=1/\uuud_{123}$. The concordance
conditions \mythetag{1.9} and \mythetag{1.10} are postulated to be 
valid for the bundles $S\uuU\!M$ and $S\uuuU\!M$ of the Standard Model.
\head
2. A little bit of the general theory.
\endhead
    Let $U\!M$ be an arbitrary $q$-dimensional complex vector bundle over
the space-time manifold $M$. A connection in $U\!M$ is a geometric structure
which is used for to differentiate tensor fields obtaining other tensor
fields from them:
$$
\hskip -2em
\nabla\!:\,U^\varepsilon_\eta\bar U^\sigma_\zeta M\to 
U^\varepsilon_\eta\bar U^\sigma_\zeta M\otimes T^*\!M.
\mytag{2.1}
$$
The operator of covariant differential \mythetag{2.1} acts upon a 
tensorial field of the type $(\varepsilon,\eta|\sigma,\zeta)$ and 
produces a tensorial field of the type $(\varepsilon,\eta|\sigma,\zeta
|0,1)$. In other words, the operator $\nabla$ adds one lower (covariant)
index to the components of a tensor. But this additional index is
associated with the tangent bundle $TM$, not with the vector bundle 
$U\!M$ itself. It is convenient to extend the domain of the operator
\mythetag{2.1} by adding more indices associated with the tangent bundle:
$$
\hskip -2em
\nabla\!:\,U^\varepsilon_\eta\bar U^\sigma_\zeta M\otimes T^m_n M\to 
U^\varepsilon_\eta\bar U^\sigma_\zeta M\otimes T^m_{n+1}M.
\mytag{2.2}
$$
The operator of covariant differential \mythetag{2.2} acts upon a 
tensorial field of the type $(\varepsilon,\eta|\sigma,\zeta|m,n)$ and 
produces a tensorial field of the type $(\varepsilon,\eta|\sigma,\zeta
|m,n+1)$. Assume that $(U,\,\boldsymbol\Psi_1,\,\ldots,\,\boldsymbol
\Psi_q)$ is a frame of the bundle $U\!M$ and let $(U,\,\boldsymbol
\Upsilon_0,\,\boldsymbol\Upsilon_1,\,\boldsymbol\Upsilon_2,\,\boldsymbol
\Upsilon_3)$ be some frame of the tangent bundle $TM$. \pagebreak Assume
that the domain $U\subset M$ is equipped with some local coordinates
$x^0,\,x^1,\,x^2,\,x^3$. Then we have the expansion
$$
\hskip -2em
\boldsymbol\Upsilon_i=\sum^3_{j=0}\Upsilon^j_i
\ \frac{\partial}{\partial x^j}
\mytag{2.3}
$$
for the frame vector fields $\boldsymbol\Upsilon_0,\,\boldsymbol\Upsilon_1,
\,\boldsymbol\Upsilon_2,\,\boldsymbol\Upsilon_3$. Apart from \mythetag{2.3},
we consider the following commutation relationships for these vector fields:
$$
\hskip -2em
[\boldsymbol\Upsilon_i,\,\boldsymbol\Upsilon_j]=
\sum^3_{k=0}c^{\,k}_{ij}\,\boldsymbol\Upsilon_k.
\mytag{2.4}
$$
The quantities $c^{\,k}_{ij}$ in \mythetag{2.4} are similar to 
structural constants of Lie algebras. For this reason they are called 
the {\it structural constants} of the frame $(U,\,\boldsymbol\Upsilon_0,
\,\boldsymbol\Upsilon_1,\,\boldsymbol\Upsilon_2,\,\boldsymbol\Upsilon_3)$,
though actually they are not constants, but smooth real-valued functions
within the domain $U$.\par
    Let $\bold X\in U^\varepsilon_\eta\bar U^\sigma_\zeta M\otimes T^m_nM$
be a tensor field of the type $(\varepsilon,\eta|\sigma,\zeta|m,n)$. In the
frame pair $(U,\,\boldsymbol\Psi_1,\,\ldots,\,\boldsymbol\Psi_q)$ and
$(U,\,\boldsymbol\Upsilon_0,\,\boldsymbol\Upsilon_1,\,\boldsymbol\Upsilon_2,
\,\boldsymbol\Upsilon_3)$ this tensor field is represented by its components
$X^{i_1\ldots\,i_\varepsilon\,\bar i_1\ldots\,\bar i_\sigma\,h_1\ldots
\,h_m}_{j_1\ldots\,j_\eta\,\bar j_1\ldots\,\bar j_\zeta\,k_1\ldots\, k_n}$.
Then $\nabla$ is represented by the formula
$$
\hskip -5em
\gathered
\nabla_{\!k_{n+1}}X^{i_1\ldots\,i_\varepsilon\,\bar i_1\ldots\,
\bar i_\sigma\,h_1\ldots\,h_m}_{j_1\ldots\,j_\eta\,\bar j_1\ldots
\,\bar j_\zeta\,k_1\ldots\, k_n}=\sum^3_{s=0}\Upsilon^s_{k_{n+1}}
\,\frac{\partial X^{i_1\ldots\,i_\varepsilon\,\bar i_1\ldots\,
\bar i_\sigma\,h_1\ldots\,h_m}_{j_1\ldots\,j_\eta\,\bar j_1\ldots
\,\bar j_\zeta\,k_1\ldots\, k_n}}{\partial x^s}\,-\\
\vspace{2ex}
\gathered
\kern -9em
+\sum^\varepsilon_{\mu=1}\sum^q_{v_\mu=1}\Alpha^{i_\mu}_{k_{n+1}\,v_\mu}\ 
X^{i_1\ldots\,v_\mu\,\ldots\,i_\varepsilon\,\bar i_1\ldots\,\bar i_\sigma
\,h_1\ldots\,h_m}_{j_1\ldots\,\ldots\,\ldots\,j_\eta\,\bar j_1\ldots\,
\bar j_\zeta\,k_1\ldots\,k_n}\,-\\
\kern 9em-\sum^\eta_{\mu=1}\sum^q_{w_\mu=1}\Alpha^{w_\mu}_{k_{n+1}\,j_\mu}
\ X^{i_1\ldots\,\ldots\,\ldots\,i_\varepsilon\,\bar i_1\ldots\,\bar i_\sigma
h_1\ldots\,h_m}_{j_1\ldots\,w_\mu\,\ldots\,j_\eta\,\bar j_1\ldots\,
\bar j_\zeta k_1\ldots\,k_n}\,+\\
\kern -9em
+\sum^\sigma_{\mu=1}\sum^q_{v_\mu=1}
\bar{\Alpha}\vphantom{\Alpha}^{\bar i_\mu}_{k_{n+1}\,v_\mu}\ 
X^{i_1\ldots\,i_\varepsilon\,\bar i_1\ldots\,v_\mu\,\ldots\,\bar i_\sigma
\,h_1\ldots\,h_m}_{j_1\ldots\,j_\eta\,\bar j_1\ldots\,\ldots\,\ldots\,
\bar j_\zeta\,k_1\ldots\,k_n}\,-\\
\kern 9em-\sum^\zeta_{\mu=1}\sum^q_{w_\mu=1}
\bar{\Alpha}\vphantom{\Alpha}^{w_\mu}_{k_{n+1}\,\bar j_\mu}\
X^{i_1\ldots\,i_\varepsilon\,\bar i_1\ldots\,\ldots\,\ldots\,\bar i_\sigma
\,h_1\ldots\,h_m}_{j_1\ldots\,j_\eta\,\bar j_1\ldots\,w_\mu\,\ldots\,
\bar j_\zeta\,k_1\ldots\,k_n}\,+\\
\kern -9em+\sum^m_{\mu=1}\sum^3_{v_\mu=0}\Gamma^{h_\mu}_{k_{n+1}\,v_\mu}\ 
X^{i_1\ldots\,i_\varepsilon\,\bar i_1\ldots\,\bar i_\sigma\,
h_1\ldots\,v_\mu\,\ldots\,h_m}_{j_1\ldots\,j_\eta\,\bar j_1\ldots\,
\bar j_\zeta\,k_1\ldots\,\ldots\,\ldots\,k_n}\,-\\
\kern 9em-\sum^n_{\mu=1}\sum^3_{w_\mu=0}\Gamma^{w_\mu}_{k_{n+1}\,k_\mu}\
X^{i_1\ldots\,i_\varepsilon\,\bar i_1\ldots\,\bar i_\sigma\,
h_1\ldots\,\ldots\,\ldots\,h_m}_{j_1\ldots\,j_\eta\,\bar j_1\ldots\,
\bar j_\zeta\,k_1\ldots\,w_\mu\,\ldots\,k_n}.
\endgathered\kern 4em
\endgathered\hskip -4em
\mytag{2.5}
$$
The quantities $\Upsilon^s_{k_{n+1}}$ for \mythetag{2.5} are taken from 
the expansion \mythetag{2.3}. In \mythetag{2.5} these quantities form
the so-called {\it Lie derivative} $L_{\boldsymbol\Upsilon_k}$:
$$
\pagebreak
\hskip -2em
L_{\boldsymbol\Upsilon_k}(f)=\sum^3_{s=0}\Upsilon^s_k
\,\frac{\partial f}{\partial x^s}.
\mytag{2.6}
$$
The quantities $\Alpha^{\,i}_{kj}$, $\bar{\Alpha}\kern -0.2pt
\vphantom{\Alpha}^{\,i}_{kj}$, and $\Gamma^{\,i}_{\!kj}$ determine 
the coordinate representation of a connection in the frame pair $(U,
\,\boldsymbol\Psi_1,\,\ldots,\,\boldsymbol\Psi_q)$ and $(U,\,
\boldsymbol\Upsilon_0,\,\boldsymbol\Upsilon_1,\,\boldsymbol\Upsilon_2,
\,\boldsymbol\Upsilon_3)$.
\mydefinition{2.1} A connection $(\Gamma,\Alpha,\bar{\Alpha})$
of the bundle $U\!M$ is called a {\it real connection\/} if the 
corresponding covariant differential \mythetag{2.2} commute with 
the involution of complex conjugation $\tau$, i\.\,e\. if $\nabla(\tau(\bold
X))=\tau(\nabla\bold X)$ for any tensor field $\bold X$.
\enddefinition
\noindent
In the case of a real connection $(\Gamma,\Alpha,\bar{\Alpha})$ we have 
the following relationships:
$$
\xalignat 2
&\hskip -2em
\Gamma^i_{kj}=\overline{\Gamma^{\kern 0.2pt\raise 0.6pt
\hbox{$\ssize i$}}_{kj}},
&&\bar{\Alpha}\vphantom{\Alpha}^i_{kj}=\overline{\Alpha^i_{kj}}.
\mytag{2.7}
\endxalignat
$$
With the use of the $\Gamma$-components of a real connection we define the
real quantities
$$
\hskip -2em
T^k_{ij}=\Gamma^k_{ij}-\Gamma^k_{j\,i}-c^{\,k}_{ij}.
\mytag{2.8}
$$
The quantities \mythetag{2.8} are the components of a tensor $\bold T$. 
It is called the {\it torsion tensor}. If $\bold T=0$, then we say that
$(\Gamma,\Alpha,\bar{\Alpha})$ is a {\it torsion-free\/} connection or a
{\it symmetric\/} connection.
\mydefinition{2.2} A real connection $(\Gamma,\Alpha,\bar{\Alpha})$ 
of the bundle $U\!M$ is called {\it concordant with the metric tensor
$\bold g$} if $\nabla\bold g=0$. If it is symmetric, i\.\,e\. if 
$\bold T=0$, then such a connection is called a {\it metric connection}.
\enddefinition
    The $\Gamma$-components of a metric connection are uniquely determined
by the metric tensor $\bold g$. In a coordinate representation we have the
following formula:
$$
\hskip -2em
\gathered
\Gamma^k_{\!ij}=\sum^3_{r=0}\frac{g^{\kern 0.5pt kr}}{2}
\left(L_{\boldsymbol\Upsilon_{\!i}}\!(g_{jr})+L_{\boldsymbol
\Upsilon_{\!j}}\!(g_{ri})-L_{\boldsymbol\Upsilon_{\!r}}\!(g_{ij})
\right)-\\
-\,\frac{c^{\,k}_{ij}}{2}
+\sum^3_{r=0}\sum^3_{s=0}g^{kr}\,\frac{c^{\,s}_{i\kern 0.5pt r}}{2}\,
g_{sj}+\sum^3_{r=0}\sum^3_{s=0}g^{kr}\,\frac{c^{\,s}_{j\kern 0.5pt r}}
{2}\,g_{s\kern 0.5pt i}.
\endgathered
\mytag{2.9}
$$
The proof of the formula \mythetag{2.9} can be found in \mycite{8},
while the formula \mythetag{2.8} for the components of the torsion 
tensor $\bold T$ is taken from \mycite{9} (see formula \thetag{6.22}
there).\par
     Each connection $(\Gamma,\Alpha,\bar{\Alpha})$ of the bundle $U\!M$
produces three curvature tensors. The components of these curvature tensors
are given by the following formulas:
$$
\align
&\hskip -5em
R^p_{kij}=L_{\boldsymbol\Upsilon_i}(\Gamma^p_{\!j\,k})
-L_{\boldsymbol\Upsilon_j}(\Gamma^p_{\!i\,k})
+\sum^3_{h=0}\left(\Gamma^p_{\!i\,h}\,\Gamma^h_{\!j\,k}
-\Gamma^p_{\!j\,h}\,\Gamma^h_{\!i\,k}\right)
-\sum^3_{h=0}c^{\,h}_{ij}\,\Gamma^p_{hk},\hskip -2em
\mytag{2.10}\\
&\hskip -5em
\goth R^p_{kij}=L_{\boldsymbol\Upsilon_i}(\Alpha^p_{j\,k})
-L_{\boldsymbol\Upsilon_j}(\Alpha^p_{i\,k})
+\sum^q_{h=1}\left(\Alpha^p_{i\,h}\,\Alpha^{\!h}_{j\,k}
-\Alpha^p_{j\,h}\,\Alpha^{\!h}_{i\,k}\right)
-\sum^3_{h=0}c^{\,h}_{ij}\,\Alpha^p_{hk},\hskip -2em
\mytag{2.11}\\
&\hskip -5em
\bar{\goth R}^p_{kij}=L_{\boldsymbol\Upsilon_i}(\bar{\Alpha}
\vphantom{\Alpha}^p_{j\,k})-L_{\boldsymbol\Upsilon_j}(\bar{\Alpha}
\vphantom{\Alpha}^p_{i\,k})
+\sum^q_{h=1}\left(\bar{\Alpha}\vphantom{\Alpha}^p_{i\,h}
\,\bar{\Alpha}\vphantom{\Alpha}^{\!h}_{j\,k}-\bar{\Alpha}
\vphantom{\Alpha}^p_{j\,h}\,\bar{\Alpha}\vphantom{\Alpha}^{\!h}_{i\,k}
\right)-\sum^3_{h=0}c^{\,h}_{ij}\,\bar{\Alpha}
\vphantom{\Alpha}^p_{hk}.\hskip -2em
\mytag{2.12}
\endalign
$$
In the case of a real connection $(\Gamma,\Alpha,\bar{\Alpha})$ the first
curvature tensor \mythetag{2.10} is a real tensor field, i\.\,e\. its 
components are real functions:
$$
\pagebreak
\hskip -2em
R^p_{kij}=\overline{R^p_{kij}}.
\mytag{2.13}
$$
Similarly, in the case of a real connection $(\Gamma,\Alpha,\bar{\Alpha})$
the third curvature tensor \mythetag{2.12} is expressed through the second 
curvature tensor \mythetag{2.11}:
$$
\hskip -2em
\bar{\goth R}^p_{kij}=\overline{\goth R^p_{kij}}.
\mytag{2.14}
$$
The formulas \mythetag{2.13} and \mythetag{2.14} are easily derived 
from \mythetag{2.7} since the parameters $c^{\,k}_{ij}$ introduced 
by the relationships \mythetag{2.4} are real-valued functions.
\head
3. The $\MatGrU(1)$-connections.
\endhead
    Having resumed some facts from the general theory, now we return
back to three special bundles considered in section~1. Let's begin with
the bundle $\uU\!M$. A metric connection for this bundle is denoted
by $(\Gamma,\uA,\bar{\uA})$. The $\Gamma$-components of this connection
are given by the formula \mythetag{2.9}.
\mydefinition{3.1} A real connection $(\Gamma,\uA,\bar{\uA})$ of the 
bundle $\uU\!M$ is called {\it concordant} with the Hermitian metric 
tensor $\bolduD$ if $\nabla\bolduD=0$.
\enddefinition
\noindent
In a coordinate form the concordance condition $\nabla\bolduD=0$ is 
written as
$$
\hskip -2em
L_{\boldsymbol\Upsilon_k}(\uD_{11})-\uD_{11}\,\uA^1_{k1}
-\uD_{11}\,\bar{\uA}\vphantom{\uA}^1_{k1}=0
\mytag{3.1}
$$
(see \mythetag{2.5}). The equality \mythetag{3.1} does not
fix the $\uA$-components of the connection $(\Gamma,\uA,\bar{\uA})$. 
It is equivalent to the following formula for their real parts:
$$
\hskip -2em
\Re(\uA^1_{k1})=\frac{L_{\boldsymbol\Upsilon_k}(\uD_{11})}
{2\,\uD_{11}}.
\mytag{3.2}
$$
Their imaginary parts are not fixed at all. Assume that $(U,
\bolduPsi_1)$ is an orthonormal frame in the sense of the
definition~\mythedefinition{1.2}. Then \mythetag{3.2} simplifies to
$\Re(\uA^1_{k1})=0$.\par
     Let $(U,\bolduPsi_1)$ and $(\tilde U,\tilde{\bolduPsi}_1)$ be two
orthonormal frames in the sense of the definition~\mythedefinition{1.2}
and assume that $U\cap\tilde U\neq\varnothing$. Then at each point
$p\in U\cap\tilde U$ we have
$$
\hskip -2em
\tilde{\bolduPsi}_1=e^{i\phi}\,\bolduPsi_1.
\mytag{3.3}
$$
If $\boldsymbol\psi$ is a wave-function being a vector with respect to the
bundle $\uU\!M$, then 
$$
\xalignat 2
&\hskip -2em
\boldsymbol\psi=\psi^1\,\bolduPsi_1,
&&\boldsymbol\psi=\tilde\psi^1\,\tilde{\bolduPsi}_1.
\mytag{3.4}
\endxalignat
$$
From \mythetag{3.3} and \mythetag{3.4} we derive the following formula for
the coefficients $\psi^1$ and $\tilde\psi^1$:
$$
\hskip -2em
\psi^1=e^{i\phi}\,\tilde\psi^1.
\mytag{3.5}
$$
Under the frame transformation \mythetag{3.3} the connection components
are transformed according to the formula similar to \thetag{19.8} in
\mycite{10}:
$$
\hskip -2em
\uA^1_{k1}=\tilde{\uA}\vphantom{\uA}^1_{k1}
-i\,L_{\boldsymbol\Upsilon_k}(\phi)
\mytag{3.6}
$$
The equality \mythetag{3.6} is an analog of the gauge transformation
for the vector-potential of the electromagnetic field (see \mycite{11}).
\par
     Now let's apply the general formula \mythetag{2.11} to the 
connection $(\Gamma,\uA,\bar{\uA})$ given by the formula \mythetag{3.3}.
As a result we get the curvature tensor with the components
$$
\hskip -2em
\goth R^1_{1ij}=
L_{\boldsymbol\Upsilon_i}(\uA^1_{j1})-L_{\boldsymbol\Upsilon_j}(\uA^1_{i1})
-\sum^3_{h=0}c^{\,h}_{ij}\,\uA^1_{h1}.
\mytag{3.7}
$$
It is the feature of one-dimensional bundles that the upper and lower 
indices do cancel each other (a contraction is performed without
summation). Therefore, the quantities \mythetag{3.7} can be understood 
as the components of a purely spacial tensor: 
$$
\hskip -2em
\goth R_{ij}=
L_{\boldsymbol\Upsilon_i}(\uA^1_{j1})-L_{\boldsymbol\Upsilon_j}(\uA^1_{i1})
-\sum^3_{h=0}c^{\,h}_{ij}\,\uA^1_{h1}.
\mytag{3.8}
$$
The quantities \mythetag{3.8} are the components of a physical field
corresponding to the gauge field with the components $\uA^1_{i1}$.
\head
4. The $\MatGrSU(2)$-connections.
\endhead
    Let $(\Gamma,\uuA,\bar{\uuA})$ be a metric connection for the bundle
$S\uuU\!M$. Its $\Gamma$-components are given by the formula \mythetag{2.9}
like in the previous case of metric $\MatGrU(1)$-connections.
\mydefinition{4.1} A real connection $(\Gamma,\uuA,\bar{\uuA})$ of the 
bundle $S\uuU\!M$ is called {\it concordant} with the Hermitian metric 
tensor $\bolduuD$ and with the skew-symmetric metric tensor $\bolduud$ 
if $\nabla\bolduuD=0$ and if $\nabla\bolduud=0$.
\enddefinition
\noindent
In a coordinate form the conditions $\nabla\bolduuD=0$ and $\nabla
\bolduud=0$ are written as
$$
\gather
\hskip -2em
L_{\boldsymbol\Upsilon_k}(\uuD_{i\bar j})-\sum^2_{a=1}\uuD_{a\bar j}
\,\uuA^a_{ki}-\sum^2_{\bar a=1}\uuD_{i\bar a}
\,\bar{\uuA}\vphantom{\uuA}^{\bar a}_{k\bar j}=0,
\mytag{4.1}\\
\hskip -2em
L_{\boldsymbol\Upsilon_k}(\uud_{ij})-\sum^2_{a=1}
\uud_{aj}\,\uuA^a_{ki}-\sum^2_{a=1}\uud_{ia}\,\uuA^a_{kj}=0
\mytag{4.2}
\endgather
$$
(see formula \mythetag{2.5}). Let's write the concordance conditions
\mythetag{4.1} and \mythetag{4.2} in some orthonormal frame $(U,\,
\bolduuPsi_1,\,\bolduuPsi_2)$ (see definition~\mythedefinition{1.3}):
$$
\xalignat 2
&\hskip -2em
\overline{\uuA^{\,i}_{\kern -0.2pt kj}}=-\uuA^{\,j}_{\kern -0.2pt k
\kern 0.2pt i},
&&\sum^2_{i=1}\uuA^{\,i}_{\kern -0.2pt k\kern 0.2pt i}=0.
\mytag{4.3}
\endxalignat
$$
The equalities \mythetag{4.3} mean that for each fixed $k$ the connection
components $\uuA^i_{kj}$ are represented by skew-Hermitian traceless
matrices. Such matrices compose the Lie algebra $\MatAlgSU(2)$ associated
with the Lie group $\MatGrSU(2)$.\par
         Let $(U,\,\bolduuPsi_1,\,\bolduuPsi_2)$ and $(\tilde U,\,
\tilde{\bolduuPsi}_1,\,\tilde{\bolduuPsi}_2)$ be two orthonormal frames 
in the sense of the definition~\mythedefinition{1.3}. Assume that $U\cap
\tilde U\neq\varnothing$. Then at each point $p\in U\cap\tilde U$ we have
$$
\hskip -2em
\tilde{\bolduuPsi}_i=\sum^2_{j=1}\goth S^j_i\,\bolduuPsi_{\!j}.
\mytag{4.4}
$$
If $\boldsymbol\psi$ is a wave-function being a vector with respect to the
bundle $S\uuU\!M$, then it can be expanded in each of the above two frames:
$$
\xalignat 2
&\hskip -2em
\boldsymbol\psi=\sum^2_{i=1}\psi^i\,\bolduuPsi_i,
&&\boldsymbol\psi=\sum^2_{i=1}\tilde\psi^i\,\tilde{\bolduuPsi}_i.
\mytag{4.5}
\endxalignat
$$
From \mythetag{4.4} and \mythetag{4.5} we derive the following formula for
the coefficients $\psi^i$ and $\tilde\psi^i$:
$$
\hskip -2em
\psi^i=\sum^2_{j=1}\goth S^i_j\ \tilde\psi^j.
\mytag{4.6}
$$
Under the frame transformation \mythetag{4.4} the connection components
$\uuA^{\,i}_{\kern -0.2pt k\kern 0.2pt j}$ are transformed according to 
the formula similar to \thetag{19.8} in \mycite{10}:
$$
\hskip -2em
\uuA^i_{kj}=\dsize\sum^2_{b=1}\sum^2_{a=1}
\goth S^i_a\,\goth T^b_j\ \tilde{\uuA}\vphantom{\uuA}^a_{k\,b}
+\vartheta^i_{kj}.
\mytag{4.7}
$$
The theta-parameters $\vartheta^i_{kj}$ in \mythetag{4.7} are defined
by the formula
$$
\hskip -2em
\vartheta^i_{kj}=\sum^2_{a=1}\goth S^i_a\,
L_{\boldsymbol\Upsilon_k}(\goth T^a_j)
=-\sum^2_{a=1}L_{\boldsymbol\Upsilon_k}
\!(\goth S^i_a)\,\goth T^a_j.\hskip -2em
\mytag{4.8}
$$
(see \thetag{9.31} in \mycite{10}). By $\goth T^i_j$ in  \mythetag{4.7}
and  \mythetag{4.8} we denote the components of the inverse matrix 
$\goth T=\goth S^{-1}$, while $L_{\boldsymbol\Upsilon_k}$ is the Lie 
derivative introduced in \mythetag{2.6}. Note that $\goth S$ and
$\goth T$ both are special unitary matrices, i\.\,e\. $\goth S\in
\MatGrSU(2)$ and $\goth T\in\MatGrSU(2)$.\par
    The curvature tensor for a $\MatGrSU(2)$-connection $(\Gamma,\uuA,
\bar{\uuA})$ is given by the general formula \mythetag{2.11} upon
substituting $\Alpha^i_{kj}=\uuA^i_{kj}$ into this formula:
$$
\goth R^p_{kij}=L_{\boldsymbol\Upsilon_i}(\uuA^p_{j\,k})
-L_{\boldsymbol\Upsilon_j}(\uuA^p_{i\,k})
+\sum^2_{h=1}\left(\uuA^p_{i\,h}\,\uuA^{\!h}_{j\,k}
-\uuA^p_{j\,h}\,\uuA^{\!h}_{i\,k}\right)
-\sum^3_{h=0}c^{\,h}_{ij}\,\uuA^p_{hk}.\quad
\mytag{4.9}
$$
The quantities \mythetag{4.9} are the components of a physical field
corresponding to the gauge field with the components $\uuA^i_{kj}$.
\head
5. The $\MatGrSU(3)$-connections.
\endhead
     The $\MatGrSU(3)$-connections correspond to the Quantum Chromodynamics
which now is a part of the Standard Model. They are associated with the
third special bundle $S\uuuU\!M$. Let $(\Gamma,\uuuA,\bar{\uuuA})$ be a
metric connection for the bundle $S\uuuU\!M$. Its $\Gamma$-components are
given by the formula \mythetag{2.9}, while $\uuuA$-components are
interpreted as gluon fields.
\mydefinition{5.1} A real connection $(\Gamma,\uuuA,\bar{\uuuA})$ of the 
bundle $S\uuuU\!M$ is called {\it concordant} with the Hermitian metric 
tensor $\bolduuuD$ and with the skew-symmetric metric tensor $\bolduuud$ 
if $\nabla\bolduuuD=0$ and if $\nabla\bolduuud=0$.
\enddefinition
\noindent
In a coordinate form the conditions $\nabla\bolduuuD=0$ and $\nabla
\bolduuud=0$ are written as
$$
\gather
\hskip -2em
L_{\boldsymbol\Upsilon_k}(\uuuD_{i\bar j})-\sum^3_{a=1}\uuuD_{a\bar j}
\,\uuuA^a_{ki}-\sum^3_{\bar a=1}\uuuD_{i\bar a}
\,\bar{\uuuA}\vphantom{\uuuA}^{\bar a}_{k\bar j}=0,
\mytag{5.1}\\
\displaybreak
\hskip -2em
L_{\boldsymbol\Upsilon_s}(\uuud_{ijk})-\sum^3_{a=1}
\uuud_{ajk}\,\uuuA^a_{si}-\sum^3_{a=1}\uuud_{iak}\,\uuuA^a_{sj}
-\sum^3_{a=1}\uuud_{ija}\,\uuuA^a_{sk}=0.
\mytag{5.2}
\endgather
$$
(see formula \mythetag{2.5}). Let's write the concordance conditions
\mythetag{5.1} and \mythetag{5.2} in some orthonormal frame $(U,\,
\bolduuuPsi_1,\,\bolduuuPsi_2,\,\bolduuuPsi_3)$ (see 
definition~\mythedefinition{1.4}):
$$
\xalignat 2
&\hskip -2em
\overline{\uuuA^{\,i}_{\kern -0.2pt kj}}=-\uuuA^{\,j}_{\kern -0.2pt k
\kern 0.2pt i},
&&\sum^3_{i=1}\uuuA^{\,i}_{\kern -0.2pt k\kern 0.2pt i}=0.
\mytag{5.3}
\endxalignat
$$
In deriving the second formula \mythetag{5.3} we used the following 
two well-known identities for a nonzero completely skew-symmetric tensor
in a three-dimensional space:
$$
\xalignat 2
&\sum^3_{k=1}\uuud_{ijk}\,\uuud^{abk}
=\delta^a_i\,\delta^b_j-\delta^a_j\,\delta^b_i,
&&\sum^3_{j=1}\sum^3_{k=1}\uuud_{ijk}\,\uuud^{ajk}=2\,\delta^a_i.
\endxalignat
$$
The equalities \mythetag{5.3} mean that for each fixed $k$ the connection
components $\uuuA^i_{kj}$ are represented by skew-Hermitian traceless
matrices. Such matrices compose the Lie algebra $\MatAlgSU(3)$ associated
with the Lie group $\MatGrSU(3)$.\par
         Let $(U,\,\bolduuuPsi_1,\,\bolduuuPsi_2,\,\bolduuuPsi_3)$ and
$(\tilde U,\,\tilde{\bolduuuPsi}_1,\,\tilde{\bolduuuPsi}_2,\,
\tilde{\bolduuuPsi}_3)$ be two orthonormal frames in the sense of the
definition~\mythedefinition{1.4}. Assume that $U\cap\tilde U\neq\varnothing$. 
Then for each point $p$ in the intersection of the
domains $U$ and $\tilde U$ we have
$$
\hskip -2em
\tilde{\bolduuuPsi}_i=\sum^3_{j=1}\goth S^j_i\,\bolduuuPsi_{\!j}.
\mytag{5.4}
$$
If $\boldsymbol\psi$ is a wave-function being a vector with respect to the
bundle $S\uuuU\!M$, then it can be expanded in each of the above two frames:
$$
\xalignat 2
&\hskip -2em
\boldsymbol\psi=\sum^3_{i=1}\psi^i\,\bolduuuPsi_i,
&&\boldsymbol\psi=\sum^3_{i=1}\tilde\psi^i\,\tilde{\bolduuuPsi}_i.
\mytag{5.5}
\endxalignat
$$
From \mythetag{5.4} and \mythetag{5.5} we derive the following formula for
the coefficients $\psi^i$ and $\tilde\psi^i$:
$$
\hskip -2em
\psi^i=\sum^3_{j=1}\goth S^i_j\ \tilde\psi^j.
\mytag{5.6}
$$
Under the frame transformation \mythetag{5.4} the connection components
$\uuuA^{\,i}_{\kern -0.2pt k\kern 0.2pt j}$ are transformed according to 
the formula similar to \mythetag{4.7}:
$$
\hskip -2em
\uuuA^i_{kj}=\dsize\sum^3_{b=1}\sum^3_{a=1}
\goth S^i_a\,\goth T^b_j\ \tilde{\uuuA}\vphantom{\uuuA}^a_{k\,b}
+\vartheta^i_{kj}.
\mytag{5.7}
$$
The theta-parameters $\vartheta^i_{kj}$ in \mythetag{5.7} are defined
by the formula similar to \mythetag{4.8}:
$$
\hskip -2em
\vartheta^i_{kj}=\sum^3_{a=1}\goth S^i_a\,
L_{\boldsymbol\Upsilon_k}(\goth T^a_j)
=-\sum^3_{a=1}L_{\boldsymbol\Upsilon_k}
\!(\goth S^i_a)\,\goth T^a_j.\hskip -2em
\mytag{5.8}
$$
By $\goth T^i_j$ in \mythetag{5.7} and  \mythetag{5.8} we again denote the
components of the inverse matrix $\goth T=\goth S^{-1}$, while 
$L_{\boldsymbol\Upsilon_k}$ is the Lie derivative (see \mythetag{2.6}).
As for $\goth S$ and $\goth T$, they both are special unitary matrices,
i\.\,e\. $\goth S\in\MatGrSU(3)$ and $\goth T\in\MatGrSU(3)$.\par
    The curvature tensor for a $\MatGrSU(3)$-connection $(\Gamma,\uuuA,
\bar{\uuuA})$ is given by the general formula \mythetag{2.11} upon
substituting $\Alpha^i_{kj}=\uuuA^i_{kj}$ into this formula:
$$
\goth R^p_{kij}=L_{\boldsymbol\Upsilon_i}(\uuuA^p_{j\,k})
-L_{\boldsymbol\Upsilon_j}(\uuuA^p_{i\,k})
+\sum^3_{h=1}\left(\uuuA^p_{i\,h}\,\uuuA^{\!h}_{j\,k}
-\uuuA^p_{j\,h}\,\uuuA^{\!h}_{i\,k}\right)
-\sum^3_{h=0}c^{\,h}_{ij}\,\uuuA^p_{hk}.\quad
\mytag{5.9}
$$
The quantities \mythetag{5.9} are the components of a physical field
corresponding to the gauge field with the components $\uuuA^i_{kj}$.
\head
6. Concluding remarks.
\endhead
     Note that the three different bundles $\uU\!M$, $S\uuU\!M$, and
$S\uuuU\!M$ are treated above in very similar ways, especially the last
two bundles. This means that crucial differences of these bundles
are hidden deeper within the Standard Model. They will be studied in
a separate paper.\par
    The formulas \mythetag{3.5} and \mythetag{3.6} yield a gauge
transformation in the case of the bundle $\uU\!M$. Note that this
gauge transformation arises as a frame transformation \mythetag{3.3}.
The same is true for the gauge transformations \mythetag{4.6},
\mythetag{4.7} and \mythetag{5.6}, \mythetag{5.7} in the case of the
bundles $S\uuU\!M$ and $S\uuuU\!M$. They are initiated by the frame
transformations \mythetag{4.4} and \mythetag{5.4} respectively.
Remember that Lorentz transformations and their spinor companions,
including $P$ and $T$-reflections, were interpreted as frame
transformations in \mycite{12}. Using the bundles $\uU\!M$, $S\uuU\!M$, 
and $S\uuuU\!M$, now we do the same for gauge transformations associated
with $\MatGrU(1)$, $\MatGrSU(2)$ and $\MatGrSU(3)$ symmetries of the
Standard Model in a gravitation field. This is the main result of this
paper.\par


\newpage
\Refs
\ref\myrefno{1}\by Kane~G.\book Modern elementary particle physics
\publ Addison-Wesley Publishing Company\yr 1987
\endref
\ref\myrefno{2}\by Aurenche~P.\paper The Standard Model of particle physics
\publ e-print \myhref{http://arXiv.org/abs/hep-ph/9712342/}{hep-ph/9712342} 
in Electronic Archive \myEarXivlink
\endref
\ref\myrefno{3}\by Hewett~J.~L.\paper The Standard Model and why we 
believe it\publ e-print 
\myhref{http://arXiv.org/abs/hep-ph/9810316/}{hep-ph/9810316} 
in Electronic Archive \myEarXivlink
\endref
\ref\myrefno{4}\by Langacker~P.\paper Structure of the Standard Model
\publ e-print \myhref{http://arXiv.org/abs/hep-ph/0304186/}{hep-ph/0304186} 
in Electronic Archive \myEarXivlink
\endref
\ref\myrefno{5}\by Willenbrock~S.\paper Symmetries of the Standard Model
\publ e-print \myhref{http://arXiv.org/abs/hep-ph/0410370/}{hep-ph/0410370}
in Electronic Archive \myEarXivlink
\endref
\ref\myrefno{6}\by Altarelli~G.\paper The Standard Model of particle
physics\publ e-print
\myhref{http://arXiv.org/abs/hep-ph/0510281/}{hep-ph/0510281} 
in Electronic Archive \myEarXivlink
\endref
\ref\myrefno{7}\by Sharipov~R.~A.\paper The electro-weak and color 
bundles for the Standard Model in a gravitation field\publ e-print 
\myhref{http://arXiv.org/abs/math/0603611/}{math.DG/0603611} 
in Electronic Archive \myEarXivlink
\endref
\ref\myrefno{8}\by Sharipov~R.~A.\paper A note on metric connections 
for chiral and Dirac spinors\publ e-print 
\myhref{http://arXiv.org/abs/math/0602359/}{math.DG}
\myhref{http://arXiv.org/abs/math/0602359/}{/0602359}
in Electronic Archive \myEarXivlink
\endref
\ref\myrefno{9}\by Sharipov~R.~A.\paper Commutation relationships and
curvature spin-tensors for extended spinor connections\publ e-print 
\myhref{http://arXiv.org/abs/math/0512396/}{math.DG/0512396}
in Electronic Archive \myEarXivlink
\endref
\ref\myrefno{10}\by Sharipov~R.~A.\paper Spinor functions of spinors
and the concept of extended spinor fields\publ e-print 
\myhref{http://arXiv.org/abs/math/0511350/}{math.DG/0511350}
in Electronic Archive \myEarXivlink
\endref
\ref\myrefno{11}\by Sharipov~R.~A.\book Classical electrodynamics and
theory of relativity\publ Bashkir State University\publaddr Ufa\yr 1997
\moreref see also
\myhref{http://arXiv.org/abs/physics/0311011}{physics/0311011}
in Electronic Archive \myEarXivlink\ and 
\myhref{http://www.geocities.com/r-sharipov/r4-b5.htm}
{r-sharipov/r4-} \myhref{http://www.geocities.com/r-sharipov/r4-b5.htm}
{b5.htm} in \myGeoCities
\endref
\ref\myrefno{12}\by Sharipov~R.~A.\paper A note on Dirac spinors 
in a non-flat space-time of general relativity\publ e-print 
\myhref{http://arXiv.org/abs/math/0601262/}{math.DG/0601262} 
in Electronic Archive \myEarXivlink
\endref
\endRefs
\enddocument
\end